\documentclass[review]{elsarticle}
\usepackage[a4paper, left=2.5cm, right=2.5cm, top=2.5cm, bottom=2.5cm]{geometry}
\usepackage{lineno,hyperref}
\usepackage{amssymb}
\usepackage{amsmath}
\usepackage{amsthm}
\usepackage{dcolumn}
\usepackage{endnotes}
\usepackage{tabularx}
\usepackage{mathtools}
\usepackage[matrix,arrow]{xy}

\journal{......}

\newtheorem{theorem}{Theorem}[section]

\newtheorem{corollary}[theorem]{Corollary}

\theoremstyle{definition}
\newtheorem{definition}[theorem]{Definition}
\newtheorem{example}[theorem]{Example}
\newtheorem{remark}[theorem]{Remark}

\newcommand{\ws}{\Sigma^*}

\newcommand\twoheaduparrow{\mathord{\rotatebox[origin=c]{90}{$\twoheadrightarrow$}}}
\newcommand\twoheaddownarrow{\mathord{\rotatebox[origin=c]{90}{$\twoheadleftarrow$}}}
\newcommand{\dda}{\twoheaddownarrow}
\newcommand{\dua}{\twoheaduparrow}
\newcommand{\ua}{\mathord{\uparrow}}
\newcommand{\da}{\mathord{\downarrow}}
\newcommand{\minn}{\mathord{\min}}
\newcommand{\maxx}{\mathord{\max}}
\newcommand{\mn}{\mathbb N}

\begin{document}

\begin{frontmatter}

\title{The set of maximal points of an $\omega$-domain need not be a $G_\delta$-set}

%




\author[1]{Gaolin Li}
\author[2,22]{Chong Shen}
\author[3]{Kaiyun Wang} \author[1]{Xiaoyong Xi\corref{cor1}}
\ead{{xixy@yctu.edu.cn}}
\cortext[cor1]{Corresponding author.}
\author[4]{Dongsheng Zhao}

\address[1]{School of Mathematics and Statistics,	Yancheng Teachers University, Yancheng 224002, Jiangsu,  China}
\address[2]{School of Science,	Beijing University of Posts and Telecommunications,  Beijing 100876,	China}
\address[22] {Key Laboratory of Mathematics and Information Networks (Beijing University of Posts and Telecommunications), Ministry of Education, China.}
\address[3]{School of Mathematics and Statistics, Shaanxi Normal University, Xi’an 710119, Shaanxi, China}
\address[4]{Mathematics and Mathematics Education, National Institute of Education,
	Nanyang Technological University,  1 Nanyang Walk, 637616, Singapore}

\begin{abstract}
A topological space  has a domain model  if it is homeomorphic to the maximal point space $\mbox{Max}(P)$ of a domain $P$.
Lawson proved that  every Polish space $X$ has an  $\omega$-domain model $P$  and for such a model $P$,  $\mbox{Max}(P)$ is a $G_{\delta}$-set of the Scott space of $P$.
Martin (2003) then asked whether it is true that for every  $\omega$-domain $Q$, $\mbox{Max}(Q)$ is $G_{\delta}$-set of the Scott space of $Q$. In this paper, we give a negative answer to Martin's long standing open problem  by constructing a counterexample.
The  counterexample here actually shows that the answer is no even  for $\omega$-algebraic domains.
 \end{abstract}


\end{frontmatter}


The Scott topology is the most useful topology in domain theory. In general, this topology is just $T_0$. However, by taking the sets $\mbox{Max}(P)$ of all maximal points of posets $P$, with the relative Scott topology, one can obtain all $T_1$ spaces (every $T_1$ space is homeomorphic to the subspace $\mbox{Max}(P)$ of some poset $P$) \cite{zhao-xi-2018}. In general, if a space $X$ is homeomorphic to $\mbox{Max}(P)$ of a poset $P$, then we say that $P$ is a poset model of $X$. The spaces with a domain model have been studied by many authors. Lawson proved that a space $X$ is a Polish space if and only if it has an $\omega$-domain model  satisfying the Lawson condition \cite{ref1}. Moreover, if $P$ is an $\omega$-domain that satisfies the Lawson condition, then the set ${\rm Max}(P)$ is a $G_{\delta}$-set of the Scott space $\Sigma P$ of $P$. As pointed out by Martin, knowing that $\mbox{Max}(P)$ is a $G_{\delta}$-subset of $\Sigma P$ is often useful in proofs: For instance, Edalat established the well-known connection between measure
theory and the probabilistic powerdomain \cite{edalat-1998} assuming a separable metric space that embedded as a $G_{\delta}$-subset of an $\omega$-domain. Martin has established several results on $G_{\delta}$-subsets  $E\subseteq \mbox{Max}(P)$,
including (i) if $P$ is an $\omega$-domain such that ${{\rm Max}(P)}$ is regular, then ${{\rm Max}(P)}$ is a $G_{\delta}$-set; (ii)  if $P$ is an $\omega$-domain that has a countable set $C\subseteq {{\rm Max}(P)}$ such that $\ua x\cap C\neq \emptyset$ for each $x\in P\backslash {{\rm Max}(P)}$, then ${{\rm Max}(P)}$ is a $G_{\delta}$-set; (iii) for any $\omega$-ideal domain $P$, ${{\rm Max}(P)}$ is a $G_{\delta}$-set.

However, as stated by Martin in  \cite[Section 8(1)]{ref2}, the answer to the following problem is not known:

\begin{itemize}
	\item Is ${{\rm Max}(P)}$ a $G_{\delta}$-set  of $\Sigma P$ for every $\omega$-domain $P$?
\end{itemize}

In this paper,   we construct  a counterexample that gives a negative answer  to the above problem.
Our example actually  shows that  even for $\omega$-algebraic domains $P$, ${{\rm Max}(P)}$ may not be a $G_{\delta}$-set.
The results here also lead to the following problem for further study on this topic: 
\begin{enumerate}[\bf Problem: ]
	\item For which $\omega$-domain $P$, is $\mbox{Max}(P)$ a $G_{\delta}$-set?
\end{enumerate}

\section{Preliminaries}

This section is devoted to a brief review of some basic concepts and notations that will be used later. For more details, we refer the readers to \cite{ref4,goubault}.

Let $P$ be a poset.
For a subset $A$ of  $P$, we shall adopt the following standard notations:
$$\ua A=\{y\in P: \exists x\in A, x\leq y\};\ \da A=\{y\in P:\exists x\in A, y\leq x\}.$$
For each $x\in X$, we simply write $\ua x$ and $\da x$ for $\ua\{x\}$ and  $\da \{x\}$, respectively. A subset $A$ of $P$ is called a \emph{lower} (resp., an \emph{upper}) \emph{set} if $A=\da A$ (resp., $A=\ua A$).
An element $x$ is \emph{maximal}  in $A\subseteq P$, if  $A\cap \ua x=\{x\}$. The set of all maximal elements of $A$ is denoted by $\maxx A$. The set of minimal elements of $A$, denoted by $\minn A$, is  defined dually.

A nonempty subset $D$ of $P$ is \emph{directed}  if every two
elements in $D$ have an upper  bound in $D$.  For  $x,y\in P$,  $x$ is \emph{way-below} $y$, denoted by $x\ll  y$, if for each directed subset
$D$ of $P$ with  $\bigvee D$ existing,
$y\leq\bigvee D$
implies $x\leq d$ for some $d\in D$.  Denote $\dua x=\{y\in P:x\ll y\}$ and $\dda x=\{y\in P:y\ll x\}$. A poset $P$ is
\emph{continuous}, if for each $x\in P$, the set $\dda x$ is directed and
$x= \bigvee\dda x$.

An element $a\in P$ is called \emph{compact}, if $a\ll a$. The set of all compact elements of $P$ will be denoted by   $K(P)$.
Then, $P$ is called \emph{algebraic} if for each $x\in P$, the set $\{a\in K(P):a\leq x\}$ is directed and
$x=\bigvee\{a\in K(P):a\leq x\}$.
A continuous (resp., algebraic) dcpo  is also called a \emph{domain} (resp., an \emph{algebraic domain}).
A subset $B\subseteq P$ is a {\em base} of  $P$ if  for each $x\in P$, $B\cap \dda x$ is directed and $\bigvee(B\cap \dda x)$=$x$.
If $P$ has a countable base, then $P$ is called an {\em $\omega$-continuous dcpo} or \emph{$\omega$-domain}.

\begin{remark}
\begin{itemize}
	\item [(1)] A dcpo is a domain if and only if it has a base.
	\item [(2)] For each base $B$ of a domain $P$, $K(P)\subseteq B$. As a consequence, every algebraic domain $P$ has the smallest base, namely $K(P)$.
\end{itemize}	
\end{remark}

A subset $U$ of  $P$ is \emph{Scott open} if
(i) $U=\mathord{\uparrow}U$ and (ii) for each directed subset $D$ of $P$ for
which $\bigvee D$ exists, $\bigvee D\in U$ implies $D\cap
U\neq\emptyset$. All Scott open subsets of $P$ form a topology on $P$,
called the \emph{Scott topology} and denoted by $\sigma(P)$. The space $\Sigma P=(P,\sigma(P))$ is called the
\emph{Scott space} of $P$.

Let $X$ be a topological space. A  subset $G \subseteq X$ is called  a \emph{$G_{\delta}$-set}, if there exists a countable family $\{U_n:n\in \mathbb N\}$ of open sets such that  $G=\bigcap_{n\in\mathbb{N}} U_{n}$.

\section{A counterexample for Martin's problem}
We shall use the following notations:
$$\begin{array}{lll}
\mathbb N&:=& \text{the set of all nonzero natural numbers};\\
\omega &:=& \text{the ordinal (also the cardinal) of natural numbers};\\
\omega_1 &:=&\text{the ordinal (also the cardinal) of real numbers};\\
\Sigma & :=&\bigcup\{\mn^{k}:k\in\mn\cup\{\omega\}\}. \text{ In other words, }  \Sigma \text{ is the set of all nonempty finite } \\&& \text{ or countably infinite sequences of } \mathbb N;\\
\Sigma^*&:=&\{{a^*}:a\in \Sigma\},\text{ which is a copy of  }\Sigma.
\end{array}$$

Each element $a\in \Sigma$ can be represented as ${\langle}a_1, a_2, a_3,\cdots, a_k{\rangle}$ if it is a finite sequence, or ${\langle}a_n{\rangle}_{n\in\mn}$ if it is an infinite sequence.
We use the notation $\ell(a)$ to denote the length of $a$, which is defined as follows:
$$
\ell(a)=\left\{\begin{array}{ll}
k, & a={\langle}a_1, a_2, a_3,\cdots, a_k{\rangle};\\
\omega, & a={\langle}a_n{\rangle}_{n\in\mn}.
\end{array}\right.
$$
We can also define $\ell(a^*):=\ell(a)$, as the lengths of $a$ and $a^*$ are clear equal.

For any $a\in \Sigma$,
we use $a_k$ (resp., $(a^*)_k$), if it exists, to denote the $k$-th element in the sequence $a$ (resp., $a^*$).
For example, if $a:=\langle{}1, 2, 3, 5, \cdots,n,n+1, \cdots\rangle$,  it follows that $\ell(a)=\omega$ and $a=\langle{a_n}\rangle_{n\in\mn}$, where  $a_n=n$ for all $n\in\mn$.

\begin{definition}
For each $a, b\in\Sigma$, $a$ is said to be a \emph{substring} of $b$, denoted by $a\sqsubseteq b$, if either of the following cases holds:
\begin{itemize}
	\item [(s1)] $\ell(a)\in\mn$, and $a_k=b_k$ for all $1\leq k\leq \ell(a)$, i.e., $a$ is a prefix of $b$;
	\item [(s2)] $\ell(a)=\omega$, and $a=b$ (in other words, $a_k=b_k$ for all $k\in\mn$).
\end{itemize}
\end{definition}

It is straightforward to verify that the relation $\sqsubseteq$ defined on $\Sigma$ is a partial order, and some of its properties  needed  for understanding the main example are listed below.

\begin{remark}\label{rem12}
	\begin{itemize}
		\item [(1)]  $(\Sigma,\sqsubseteq)$ is an algebraic domain.
		
		In fact, for each $k\in\mathbb N$,
	the set  $\Sigma_k:=\{a\in\Sigma: a_1=k\}$
together with the substring order $\sqsubseteq$ form a tree-like dcpo, as shown in Figure \ref{fig:example1}.
		\begin{figure}[h]
			\centering
			\includegraphics[width=1\textwidth]{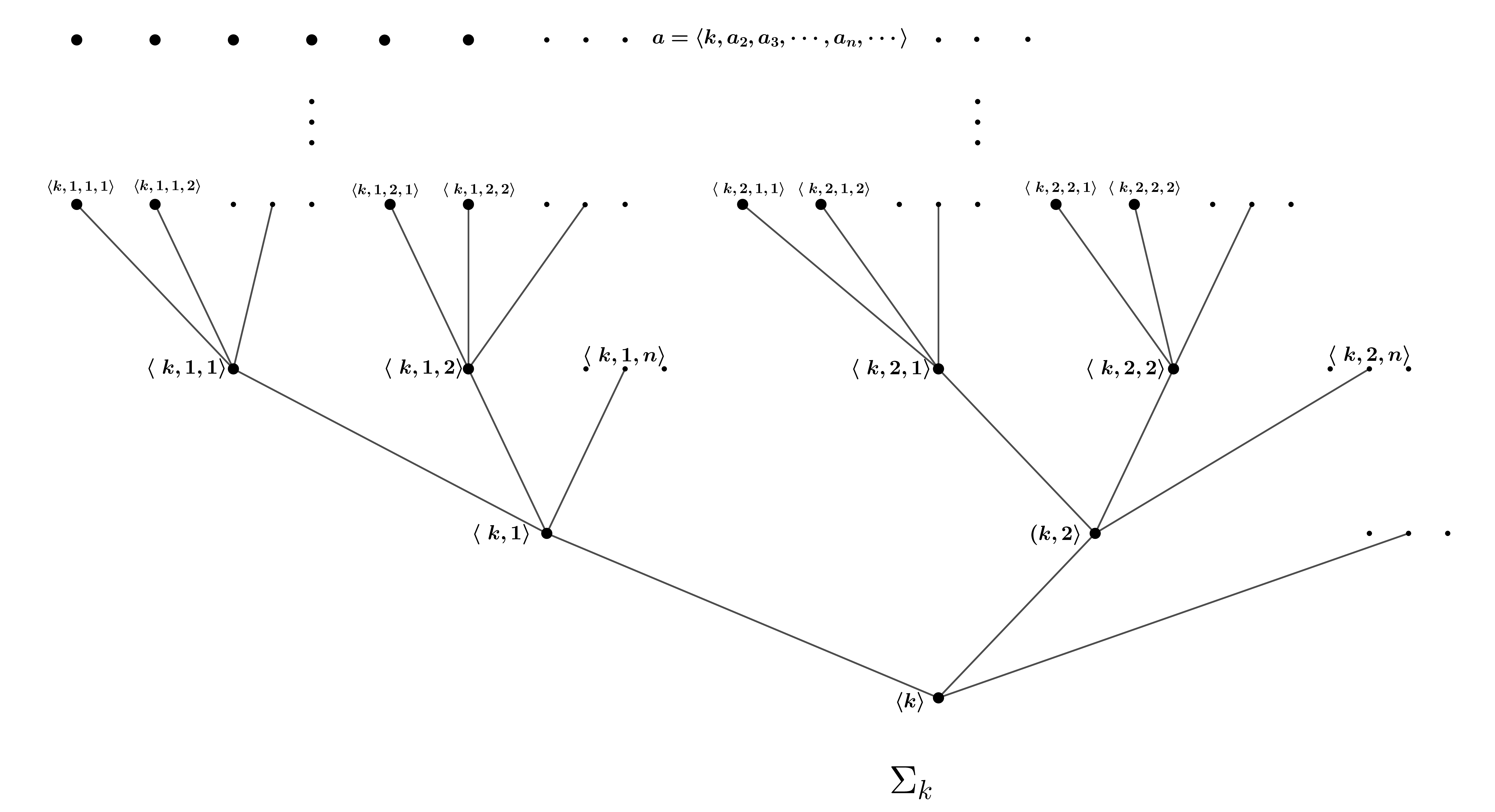}
			\caption{The subdcpo $\Sigma_k$}
			\label{fig:example1}
		\end{figure}
In addition, the poset $(\Sigma,\sqsubseteq)$ is  the countable disjoint union of these $\Sigma_k$, where $k\in\mathbb N$. That is, for any $x,y\in \Sigma$, $x\sqsubseteq y$ in $\Sigma$ iff there exists some $k_0\in \mn$ such that $x,y\in\Sigma_{k_0}$ and $x\sqsubseteq y$ in $\Sigma_{k_0}$,  as shown in Figure \ref{fig:example2}.
	\begin{figure}[h]
			\centering
			\includegraphics[width=1\textwidth]{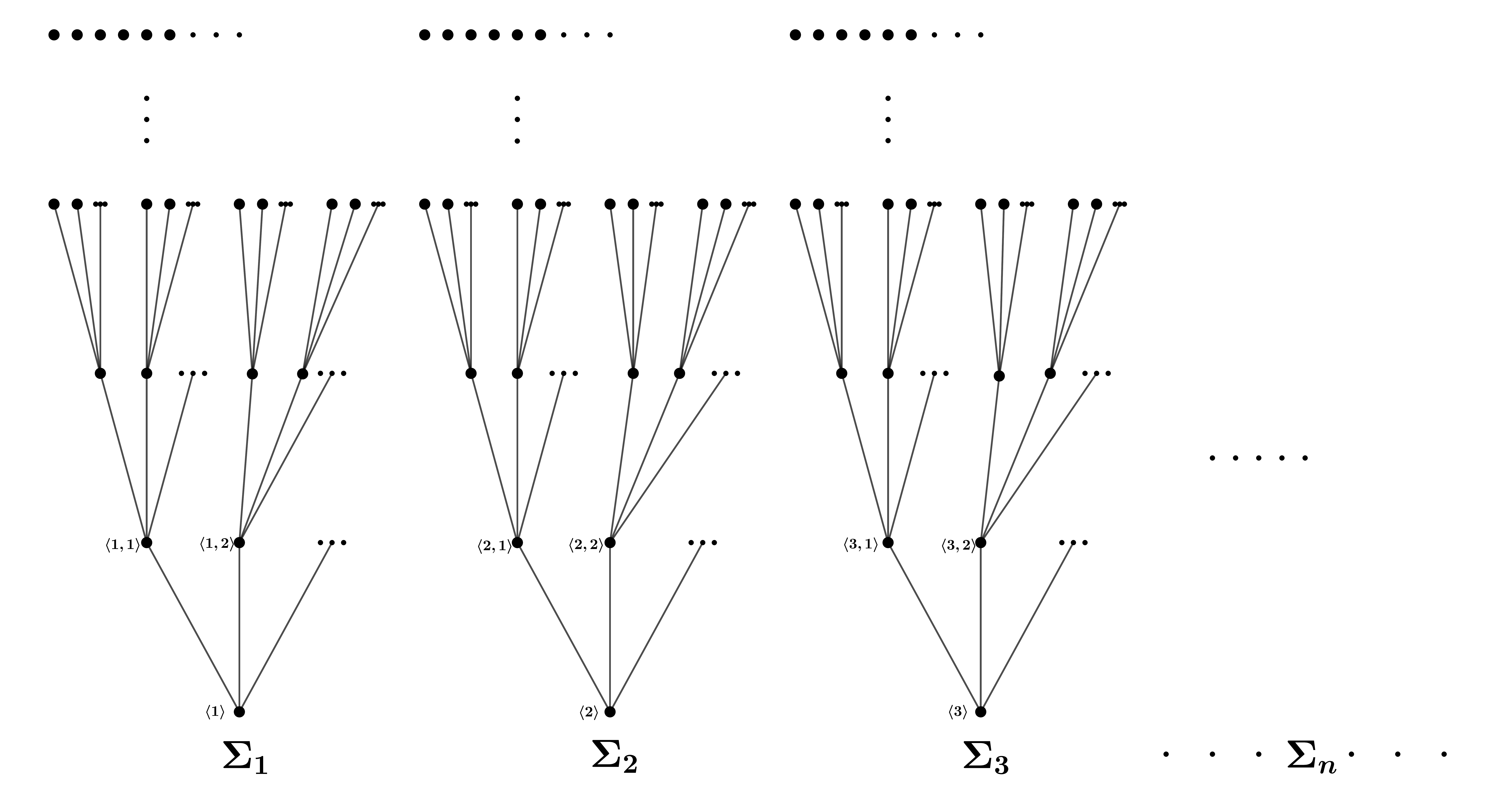}
			\caption{The algebraic domain $(\Sigma,\sqsubseteq)$}
			\label{fig:example2}
		\end{figure}
\item [(2)]  $\maxx\Sigma=\{a\in\Sigma:\ell(a)=\omega\}$, $|\Sigma|=|\maxx\Sigma|=\omega_1$
and $|\Sigma\setminus\maxx\Sigma|=\omega$.
	\item [(3)]  The compact elements of $\Sigma$ are the nonempty finite sequences of $\mn$. In other words,
		$$K(\Sigma) = \{a\in\Sigma: \ell(a)\in\mn\}.$$
		\item [(4)] For any two  elements $a,b\in \Sigma$, $\{a,b\}$ has upper bounds iff 
		$a$ and $b$ are comparable (i.e., either $a\leq b$ or $b\leq a$).
	\end{itemize}
 \end{remark}

We can also define a partial order $\sqsubseteq^*$ on $\Sigma^*$ as follows: $\forall a,b\in\Sigma$,
$$a^*\sqsubseteq b^*\mbox{ iff } a\sqsubseteq b.$$
It is clear that the posets $(\Sigma,\sqsubseteq)$ and $(\Sigma^*,\sqsubseteq^*)$ are order-isomorphic. Furthermore, the properties of $(\Sigma,\sqsubseteq)$   listed in Remark \ref{rem12} also hold for $(\Sigma^*,\sqsubseteq^*)$.

\begin{example}\label{exam1}
Let $X:=\{{x}_{m,n}:m\in  \mathbb{N},n\in \mathbb N\cup\{\omega\}\}$ with the  partial order $\leq_X$ defined as follows: $\forall m_1,m_2\in \mn$ and $n_1,n_2\in \mathbb N\cup\{\omega\}$,
$$x_{m_1,n_1}\leq_X x_{m_2,n_2} \Leftrightarrow m_1=m_2\text{ and }n_1\leq n_2.$$
Then,  $X$ is an $\omega$-algebraic domain, as shown in Figure \ref{fig:example3}.
\begin{figure}[h]
	\centering
	\includegraphics[width=0.7\textwidth]{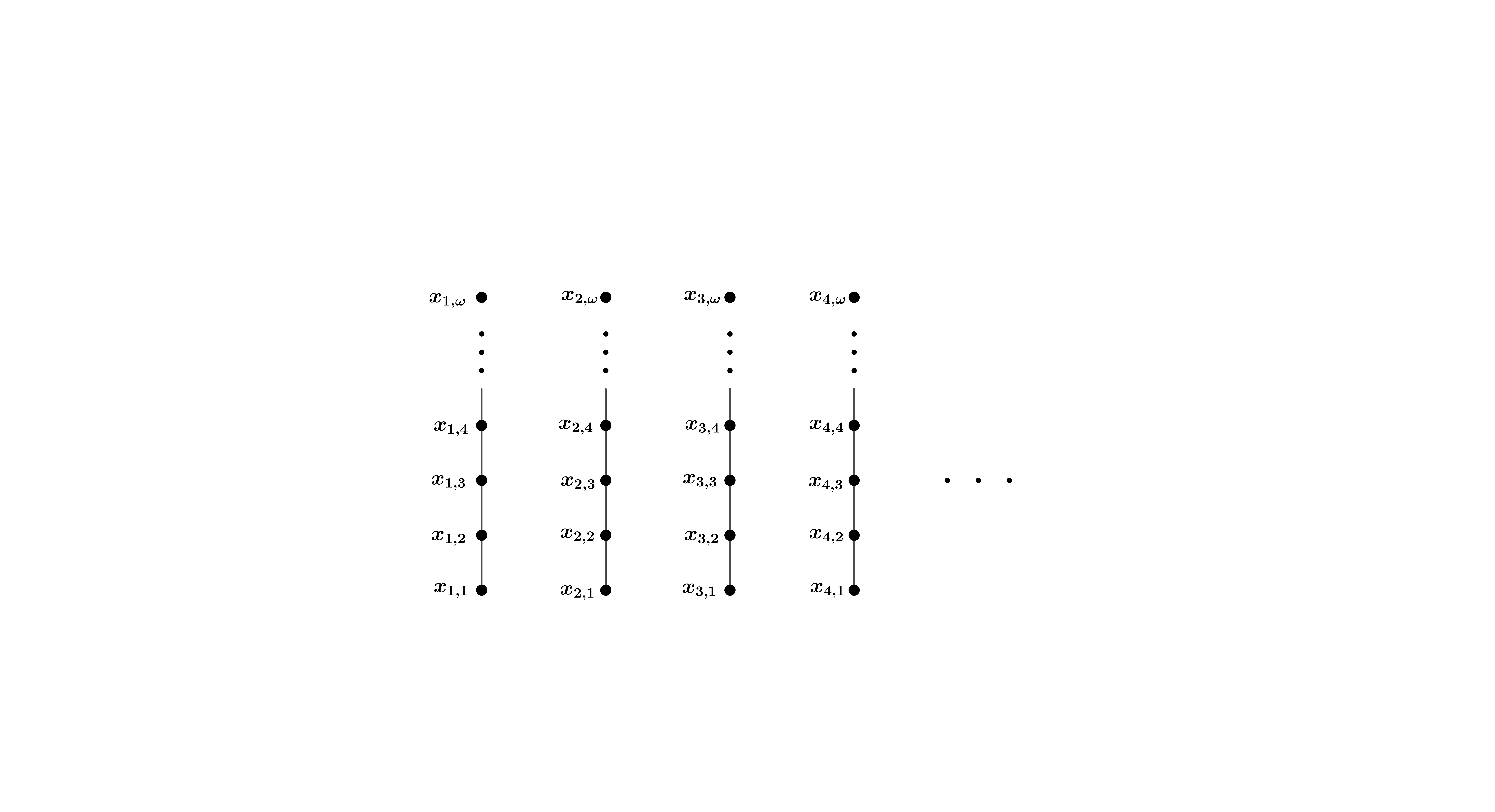}
	\caption{The algebraic domain $(X,\leq_X)$}
	\label{fig:example3}
\end{figure}
In fact, one  easily sees that $\bigvee_{n\in\mn}x_{k,n}=x_{k,\omega}$ for all $n\in\mn$. Moreover, we have that
\begin{center}
$\maxx X=\{{x}_{m,\omega}:m\in  \mathbb{N}\}$ and	$K(X)=X\setminus\maxx X=\{{x}_{m, n}:m, n\in  \mathbb{N}\}$.
\end{center}
It then follows that
$|K(X)|=|\maxx X|=|X|=\omega$.
\end{example}

Before presenting the counterexample, let us briefly discuss the order structure defined by combining two ordered structures, which will be helpful to understand the order relationship of our main counterexample.

\begin{remark}\label{rem}
Let $P=\{x_n:n\in\mn\cup\{\omega\}\}\cup\{y_n:n\in\mn\cup\{\omega\}\}$. 
We define two relations $\leq_1$ and $\leq_2$ on $P$.
\begin{itemize}
	\item [(1)] The relation $\leq_1 $ is given as follows:
	\begin{itemize}
		\item [(P1)] $x_1\leq_1 x_2\leq_1 x_3\leq_1 \cdots\leq_1 x_{n}\leq_1 x_{n+1}\leq_1\cdots\leq_1 x_{\omega}$;
		\item [(P2)] $y_1\leq_1 y_2\leq_1 y_3\leq_1 \cdots\leq_1 y_{n}\leq_1 y_{n+1}\leq_1\cdots\leq_1 y_{\omega}$;
		\item [(P3)] $\forall m,n\in\mn$, $x_{m}\leq_1 y_{n}$ $\Leftrightarrow$ $m\leq n$.
	\end{itemize}
Then, we define  $\leq_2:=\leq_1\cup\{(x_{\omega}, y_{\omega})\}$.
It is east to verify that both  $\leq_1$ and $\leq_2$ are partial orders on $P$, as shown in Figure \ref{fig:example6}.
	 \begin{figure}[h]
	\centering
	\includegraphics[width=0.5\textwidth]{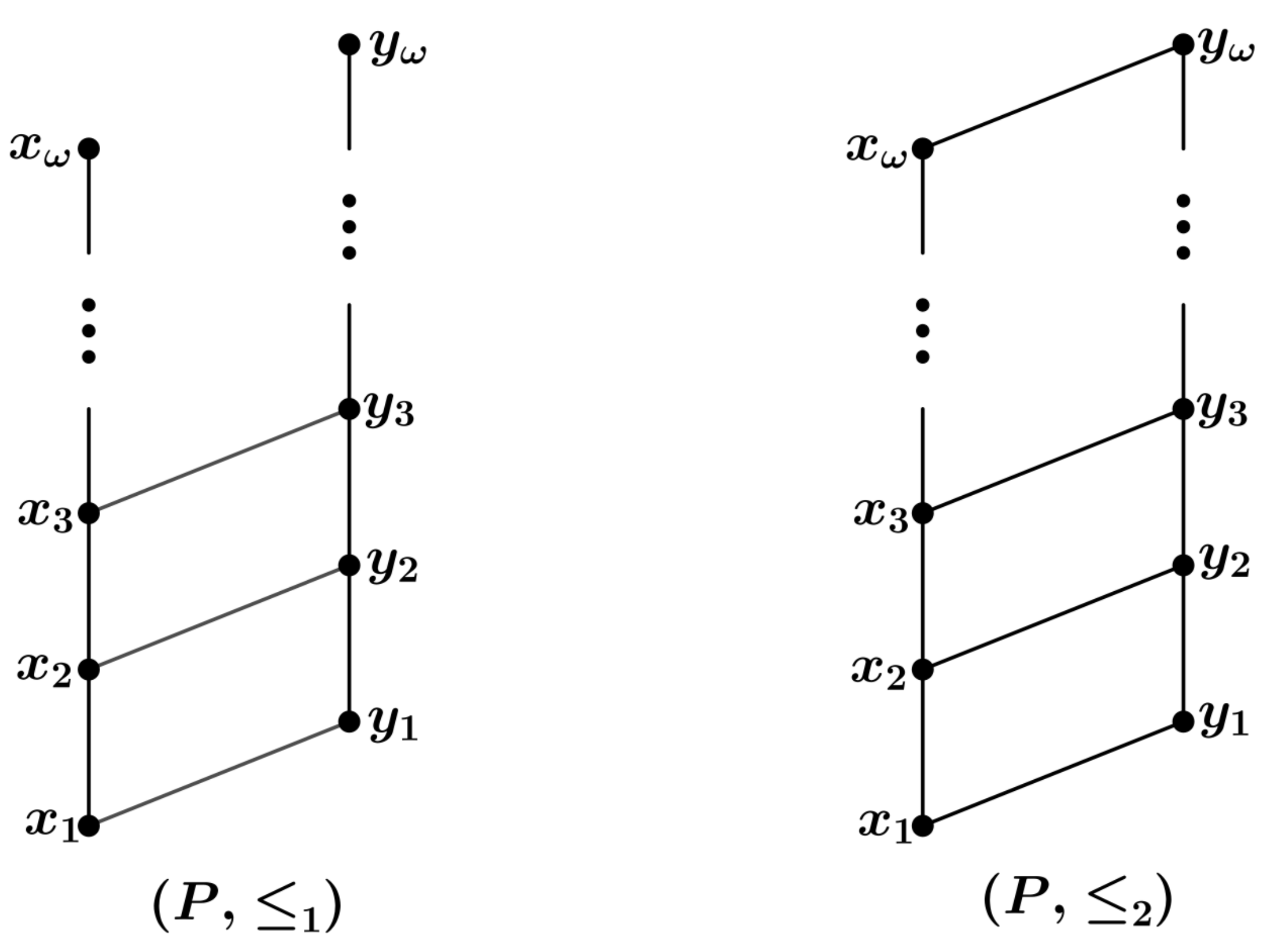}
	\caption{The posets $(P,\leq_1)$ and $(P,\leq_2)$}
	\label{fig:example6}
\end{figure}
\item [(2)]  In the poset $(P,\leq_1)$, the chain $C:=\{x_n:n\in\mn\}$ has two incomparable upper bounds, namely $x_{\omega}$ and $y_\omega$. Hence, the supremum of $C$  in $(P,\leq_1)$ does not exist, and consequently $\bigvee_{(P,\leq_1)}C\neq x_{\omega}$.
However, since $x_{\omega}\leq_2y_{\omega}$, it is clear that the supremum of $C$  in $(P,\leq_2)$ exists and equals  $x_{\omega}$, i.e.  $\bigvee_{(P,\leq_2)} C=x_{\omega}$.
\end{itemize}
\end{remark}

\begin{example}
Let
$L:=X\cup\Sigma\cup\Sigma^*,$
where $X$ is the poset defined in Example \ref{exam1}.
For any $u,v\in L$, define $u\leq v$ if either  of the following conditions holds:
\begin{itemize}
	\item [(i)] $u\sqsubseteq v$ in $\Sigma$;
	\item [(ii)]  $u\sqsubseteq^* v$ in $\Sigma^*$;
	\item [(iii)] $u\leq_X v$ in $X$;
	\item [(iv)] $u\in\Sigma$, $v=b^*\in\Sigma^*$ and $u \sqsubseteq b$, as shown in Figure \ref{fig:example44}.
	
	\medskip
	
		 \begin{figure}[h]
		\centering
		\includegraphics[width=0.7\textwidth]{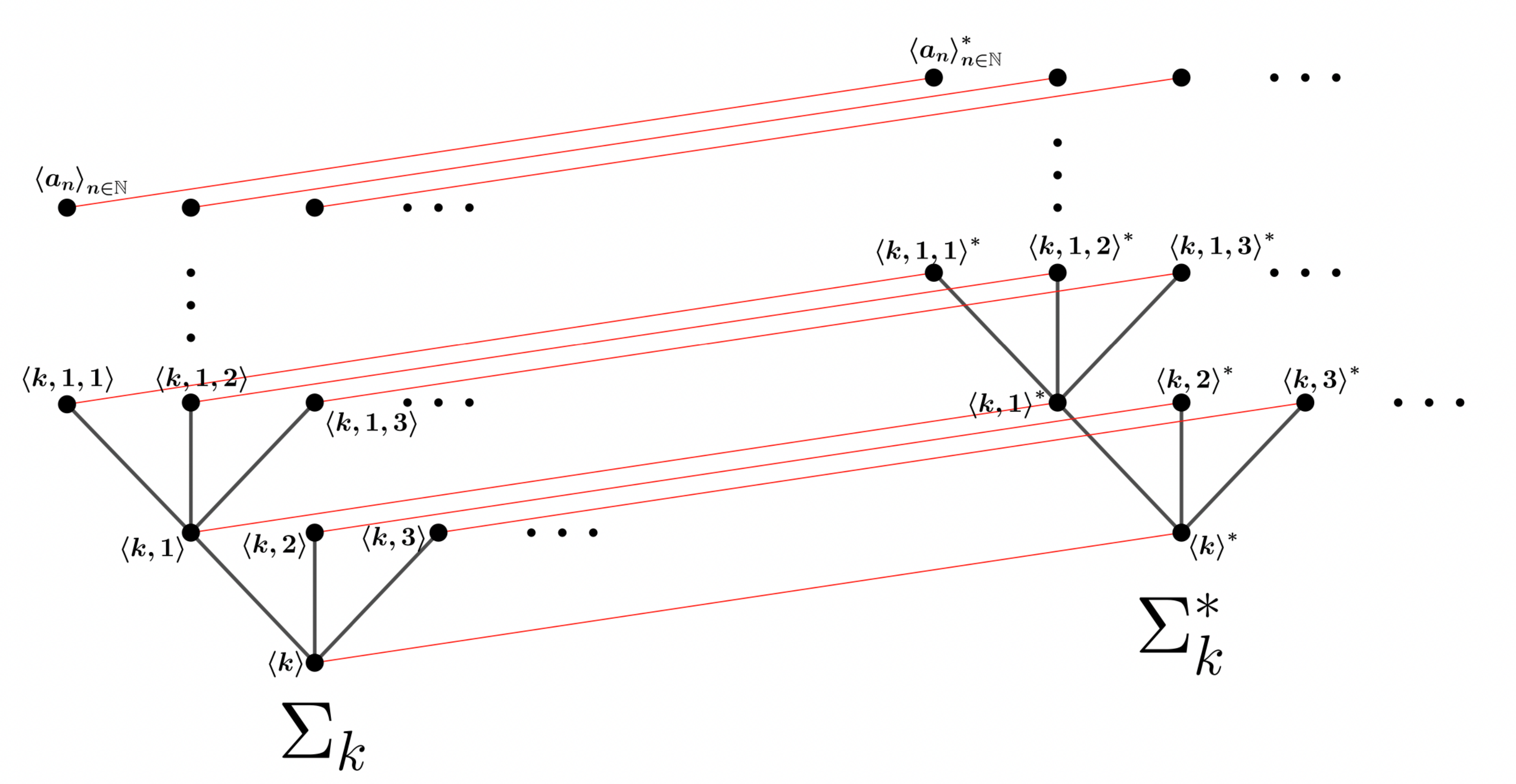}
		\caption{The relations between $\Sigma$ and $\Sigma^*$}
		\label{fig:example44}
	\end{figure}

	\item [(v)] $\exists k,m\in\mn$, $u=x_{k,m}$, $\ell(v)\geq k$ and $v_k\geq m$
	 (refer to Figure \ref{fig:example4}, where the black lines in the figure represent the relationship between the points in $X$ and $\Sigma$, while the dotted red lines represent  the relations between the points in $\Sigma$).
	 \begin{figure}[h]
	 	\centering
	 	\includegraphics[width=1\textwidth]{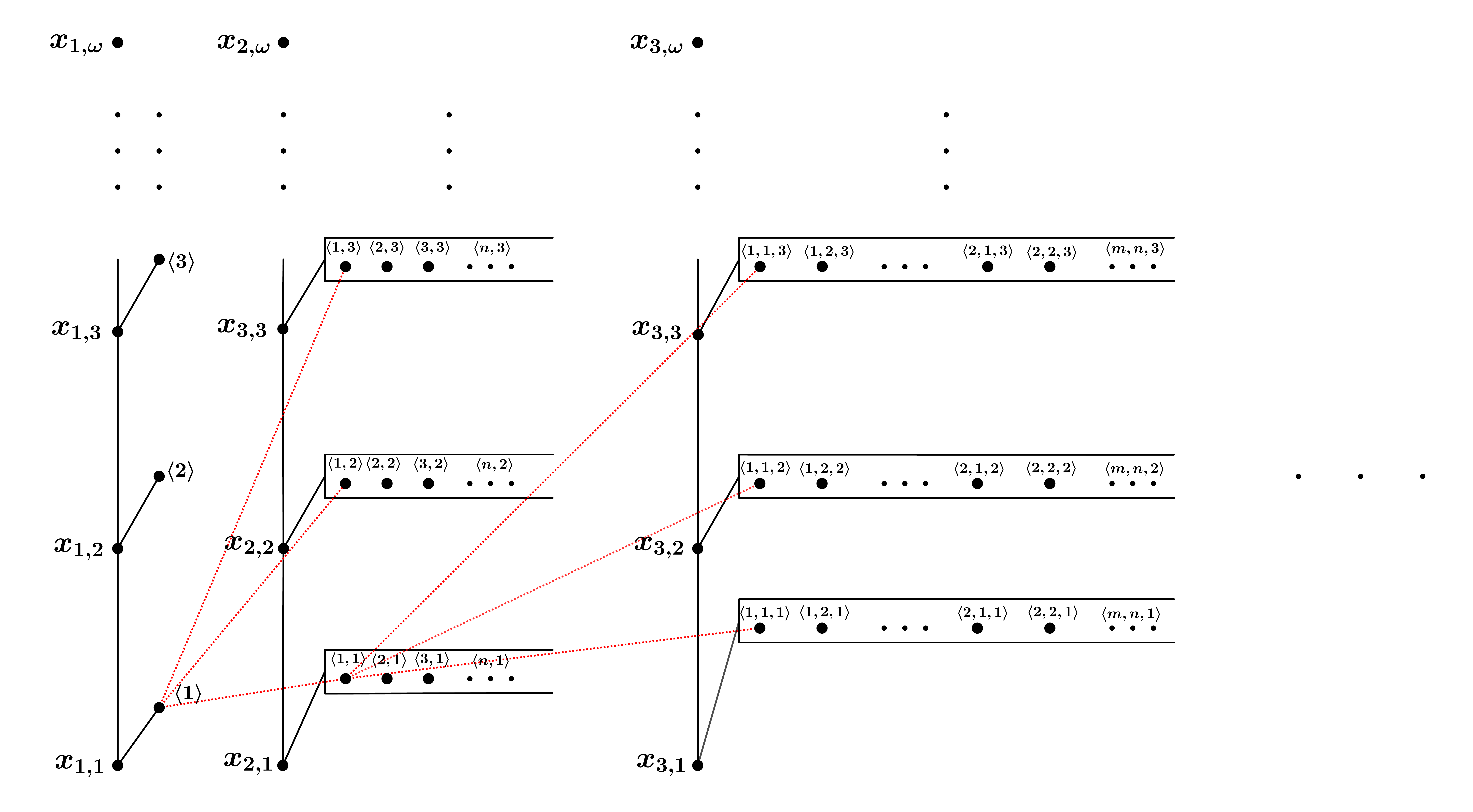}
	 	\caption{The subset $X\cup\Sigma$ of $L$}
	 		\label{fig:example4}
	 \end{figure}
\end{itemize}

In fact, Condition (v) describes the relations between $X$ and $\Sigma\cup\Sigma^*$, which can be characterized as follows: for any $k, m\in \mn$, $v\in \Sigma\cup\Sigma^*$,
\begin{center}
	$x_{k,m}\leq v$\  iff \
	$\exists n_1,n_2,\cdots,n_{k}\in\mn$ such that $n_k\geq m$ and
$x_{k,m}\leq \langle{n_1,n_2,\cdots,n_k} \rangle\leq v.$
\end{center}
	 For example, we have the following facts:
	 $$
	 \begin{array}{ll}
&x_{4,11}\leq \langle 1,5,7,11 \rangle\leq \langle
1,5,7,11,11 \rangle^*\leq\langle 1,5,7,11,11,22,22\rangle^*;\\
&x_{4,11}\leq \langle{1,5,7,11}\rangle\leq\langle 1,5,7,111 \rangle\leq \langle
1,5,7,111,11 \rangle^*\leq\langle 1,5,7,111,11,22,22\rangle^*;\\
&x_{3,3}\leq \langle{}1,2,3\rangle\leq \langle{}1,2,3,\cdots,n,n+1,\cdots\rangle\leq \langle{}1,2,3,\cdots,n,n+1,\cdots\rangle^*;\\
&x_{3,3}\leq \langle{}2,3,4\rangle\leq \langle{}2,3,4,\cdots,n,n+1,\cdots\rangle\leq \langle{}2,3,4,\cdots,n,n+1,\cdots\rangle^*.
	 \end{array}$$
	
Next, we  prove that $L$ is an $\omega$-algebraic dcpo such that $\maxx L$ is not a $G_{\delta}$-set in the Scott space of $L$. We achieve  this in a few steps.

\medskip

(1)  The relation $\leq$ is a partial order on $L$.

The reflexivity and antisymmetry are almost trivial. For showing the transitivity, we need the following:

\medskip

\emph{Claim 1: } $\forall k, m\in \mn$, $\forall a\in \Sigma$, \ $x_{k,m}\leq a$ \ $\Leftrightarrow$ \  $x_{k,m}\leq a^*$.

By the definition of $a^*$, it holds that  $\ell(a^*)=\ell(a)$ and $(a^*)_k=a_k$.
Therefore, Claim 1 follows from the following facts:
\begin{center}
	$x_{k, m}\leq a$\ iff\ $\ell(a^*)=\ell(a)\geq k$ and $(a^*)_k=a_k\geq m$
\	iff \ $x_{k, m}\leq a^*$.
\end{center}

\medskip

\emph{Claim 2: } $\forall a, b\in \Sigma$, \ $a\leq b$ \ $\Leftrightarrow$ \  $a\leq b^*$  \ $\Leftrightarrow$\ $a^*\leq b^*$\ $\Leftrightarrow$\ $a\sqsubseteq b$.

Claim 2  follows from the following facts:
$$
	a\leq b\ \stackrel{(i)}{\Leftrightarrow}\  a\sqsubseteq b \text{ in }\Sigma
\ \stackrel{(iv)}{\Leftrightarrow} \ a\leq b^*\
\stackrel{(iv)}{\Leftrightarrow}\  a\sqsubseteq b \text{ in }\Sigma\
\Leftrightarrow\  a^*\sqsubseteq^* b^* \text{ in }\Sigma^*\ \stackrel{(ii)}{\Leftrightarrow}\ a^*\leq b^*.
$$

\medskip

We now  verify the transitivity
 only by considering the following cases:
for any $k,m,k',m'\in \mn$ and $a,b,c\in \Sigma$,
\begin{itemize}
	\item [(c1)] $x_{k,m}\leq x_{k',m'}\leq a\leq b$. Then, as $x_{k,m}\leq x_{k',m'}$, we have that $k=k'$ and $m\leq m'$. Also, since $x_{k',m'}\leq a\leq b$, we have that   $\ell(b)\geq \ell(a)\geq k'=k$ and $b_{k}=b_{k'}=a_{k}=a_{k'}\geq m'\geq m$, and by (v), it follows that  $x_{k,m}\leq b$;
		\item [(c2)] $x_{k,m}\leq a\leq b$. Then, $\ell(b)\geq \ell(a)\geq k$ and $b_k=a_k\geq m$, which implies that $x_{k,m}\leq b$;
	\item [(c3)]  By Claims 1 and 2, both cases $x_{k,m}\leq a\leq b^*$ and $x_{k,m}\leq a^*\leq b^*$ are equivalent to
$x_{k,m}\leq a\leq b$. Then by (c2) , one can imply that $x_{k,m}\leq b$, and this is equivalent to $x_{k,m}\leq b^*$ by Claim 1.
\item [(c4)] By Claim 2, both cases $a\leq b^*\leq c^*$ and $a\leq b\leq c^*$ are equivalent to $a\leq b\leq c$, which means  $a\sqsubseteq b\sqsubseteq c$ in $\Sigma$. This  implies $a\sqsubseteq c$, which is equivalent to $a\leq c^*$.
\end{itemize}

\medskip

(2)  $\maxx L=\maxx \Sigma^*\cup \maxx X$ and $|\maxx L|=\omega_1$.

First, one easily observes that
$$\begin{array}{lll}
	\maxx \Sigma &=&\mathbb N^{\omega};\\
	\maxx \Sigma^* &=&\{{a^*}:a\in \mathbb N^{\omega}\};\\
	\maxx X&=&\{{x}_{m,\omega}:m\in\mn\},
\end{array}$$
and $$|\maxx \Sigma |=|\maxx \Sigma^* |=\omega_1,\ |\maxx X|=\omega.$$
In addition, by the definition of the order $\leq $ on $L$, we have that $\ua x_{m,\omega}\cap L=x_{m,\omega}$, and therefore, $\maxx X\subseteq \maxx L$. Then,   we conclude that $\maxx L=\maxx \Sigma^*\cup \maxx X$, and  $|\maxx L|=\omega_1$.

\medskip

In view of Remark \ref{rem}, we need to check the following fact.

(3)  $L$ is a dcpo.

\medskip

Before proving this, we need the following  remark about the characteristics of the relations of points $X\cup\Sigma$.

\emph{Remark: } For any  $m\in\mn$,  the set $A_m$ of all elements $a\in \Sigma$ that satisfy $x_{m,k}\leq a$ and $\ell(a)=m$ for some $k\in\mn$,  is an antichain.
That is, $A_m=\mn^m=\bigcup_{k\in\mn}\min(\ua x_{m,k}\cap\Sigma)$, is an antichain. For example, $A_3$ is the antichain of all black points (on the right side) in  Figure \ref{fig:example7}.
\begin{figure}[h]
	\centering
	\includegraphics[width=1\textwidth]{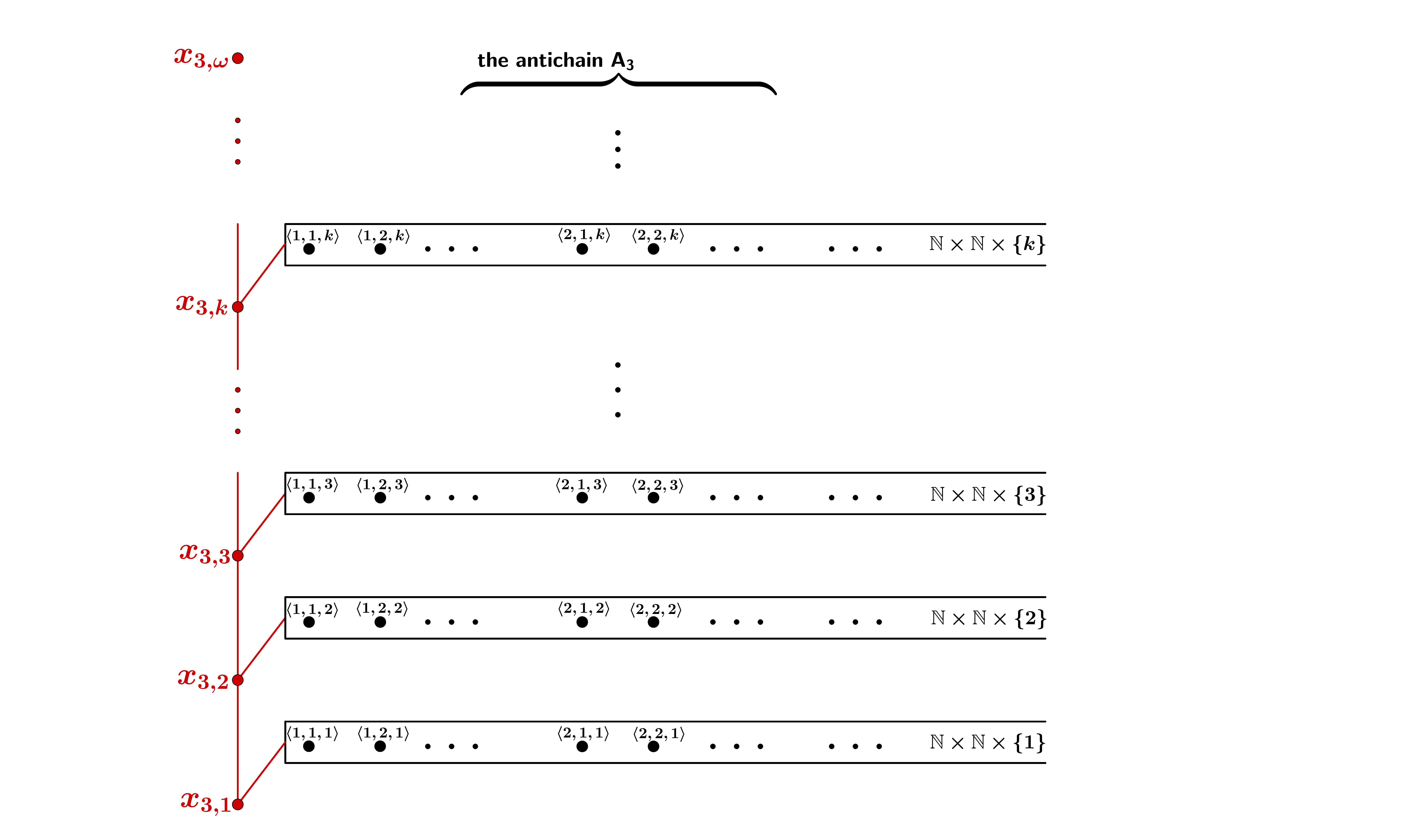}
	\caption{The set $\{x_{3,k}:k\in \mn\cup\{\omega\}\}\cup A_3$}
	\label{fig:example7}
\end{figure}
By Remark \ref{rem12}(4), we see that any two distinct elements of $A_m$ have no upper bound in $L$.
Thus,
for any  $a\in \Sigma$ and $k, l\in\mn$, we have that
$a\geq x_{m,k}$ iff there exists a unique element $b\in A_m$ such that $a\geq b$. This fact plays a crucial role in guaranteeing that the supremum of the chain $\{x_{m,k}:k\in\omega\}$ exists and is equal to $x_{m,\omega}$ (comparing to Remark \ref{rem}). The convergence of all such chains in $X$ to their respective maximal points is the key to proving that $L$ is a dcpo, as shown below in Claim 1.

\medskip

Now  $L$ being a dcpo  can be confirmed  by the following four claims.

\medskip

\emph{Claim 3:\ } For any $m\in\mn$ and any infinite subset $\{n_k:k\in\mn\}$ of $\mn$, $\bigvee_{k\in\mn}x_{m, n_k}=x_{m,\omega}$.

Comparing to Remark \ref{rem}, it is not obvious.

First, note that $x_{m,\omega}$ is an upper bound of $\{x_{m,n_k}:k\in\mn\}$ in $L$.
It remains to verify that $x_{m,\omega}$ is the unique upper bound.
We show this by considering the two cases below:

(i)  Assume $\{x_{m,n_k}:k\in\mn\}$ has an upper bound in $\Sigma$.  Then, as $\Sigma\subseteq \da \max\Sigma$, there must exist an element $a:=\langle{a_k}\rangle_{k\in\mn}\in \max\Sigma$ such that $x_{m,n_{k}}\leq a$ for all $k\in\mn$. On the one hand, as $\{n_k:k\in\mn\}$ is  infinite,
it follows that $\bigvee_{k\in\mn}n_k=\omega> a_m\in \mn$, so there exists $k_0\in \mn$ such that $n_{k_0}>a_m$. On the other hand, as
$x_{m,n_{k_0}}\leq a$, we have that $a_m\geq n_{k_0}>a_m$,  a contradiction, as shown in Figure \ref{fig:example5}.
\begin{figure}[h]
	\centering
	\includegraphics[width=0.6\textwidth]{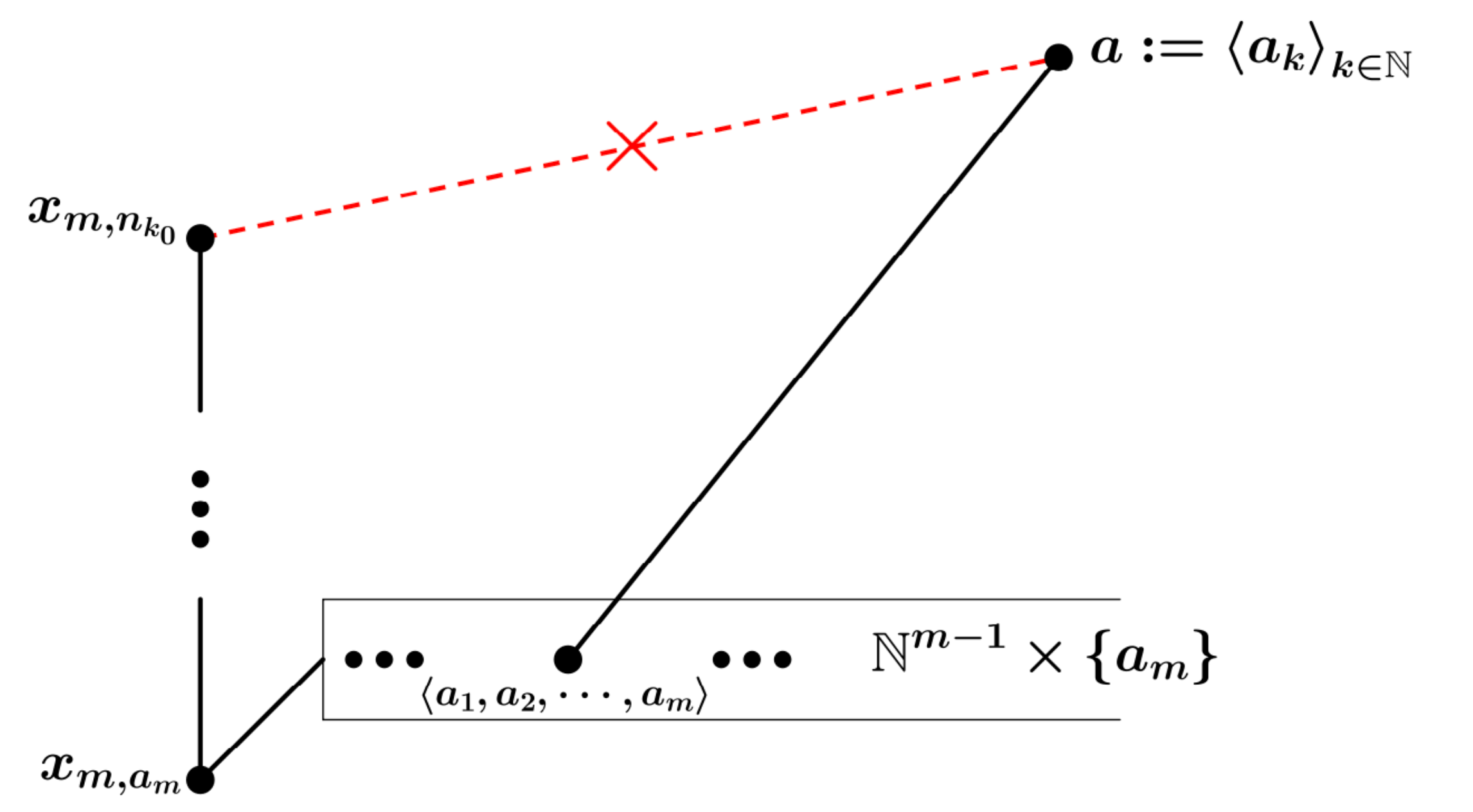}
	\caption{Order relations of $x_{m,a_m}$, $x_{m,n_{k_0}}$ and $a$}
	\label{fig:example5}
\end{figure}
Therefore, $\{x_{m,n_k}:k\in\mn\}$ has no upper bound in $\Sigma$.

(ii) Assume that  there exists $a^*\in \Sigma^*$ that is an upper bound of $\{x_{m,n_k}:k\in\mn\}$. Then, by Claim 1 in (1),  $a\in \Sigma$ is also an upper bound of  $\{x_{m,n_k}:k\in\mn\}$, which  contradicts to  (i).

Then (i) and (ii) together show that $\{x_{m,n_k}:k\in\mn\}$ has no upper bound in $\Sigma\cup\Sigma^*$, and as it is clear that the only upper bound in $X$ is $x_{m,\omega}$, we deduce that $x_{m,\omega}$ is the unique upper bound of $\{x_{m, n_k}:k\in\mn\}$  in $L$. Therefore, $\bigvee_{k\in\mn}x_{m,n_k}=x_{m,\omega}$.

\medskip

\emph{Claim 4:\ } All sets $X, \Sigma$ and $\Sigma^*$ are subdcpos of $L$.	

Clearly both $\Sigma$ and $\Sigma^*$ are subdcpos of $L$.
We now check that $X$ is a subdcpos of $L$. Note that $X$ with the restriction order $\leq$ is exactly $\leq_X$ defined in Example \ref{exam1}.
Note that all infinite directed subsets of $X$ are of the form $\{x_{m,n_k}:k\in\mn\}$, where $m\in \mn$ and $\{n_k:k\in\mn\}$ is an infinite subset of $\mn$, hence  we
only need to prove that $\bigvee_{k\in\mn}x_{m,n_k}=x_{m,\omega}$.

\medskip

\emph{Claim 5:\ } For any directed subset $D$ of $L$,  if $\bigvee D\in{\Sigma^*}$, then $D\cap \ws$ is a directed set such that $\bigvee D=\bigvee (D\cap \ws)$.

It suffices to prove $D\subseteq\da (D\cap \Sigma^*)$. Note that $X\cup \Sigma$ is both a lower set and a subdcpo of $L$, thus Scott closed.
Then,  $\Sigma^*=L\setminus(X\cup\Sigma)$ is a Scott open subset of $L$, so there exists $d_0\in D\cap \Sigma^*$, it follows that  $\ua d_0\subseteq \ua \Sigma^*=\Sigma^*$. Hence $D\subseteq \da (D\cap \ua d_0)\subseteq \da (D\cap \Sigma^*)$. This shows that $D\cap \Sigma^*$ is cofinal in $D$, and so $D\cap \Sigma^*$ is a directed set whose supremum equals $\bigvee D$ (see  pp.61 in \cite{goubault}).

\medskip

\emph{Claim 6:\ } For any directed subset $D$ of $L$,  if $\bigvee D\in\Sigma$, then $D\cap \Sigma$ is a directed set such that $\bigvee D=\bigvee (D\cap \Sigma)$.

For this, we only need to prove $D\subseteq\da (D\cap \Sigma)$. By the definition of  $\leq $ on $L$, we have that
$D\subseteq \da \bigvee D\subseteq \da \Sigma=\Sigma\cup (X\setminus\maxx X)$, and hence  $D=(D\cap \Sigma)\cup(D\cap (X\setminus\maxx X))$.
Thus, to show $D\subseteq\da (D\cap \Sigma)$, it suffices to prove
 $D\cap (X\setminus\maxx X)\subseteq \da (D\cap\Sigma)$.
  Let ${x}_{m,n}\in D\cap (X\setminus\maxx X)$.
Note that $D\cap \Sigma\neq\emptyset$ (if $D\cap \Sigma=\emptyset$, then $D\subseteq X$, and since $X$ is a subdcpo of $L$, $\bigvee D\in X$, a contradiction), so we choose an arbitrary point $b\in D\cap \Sigma$. Since $D$ is directed, there exists $e\in D$ such that $b\leq e$ and ${x}_{m,n}\leq e$. Recall that $D\subseteq \da \Sigma$ and $\ua b\cap \da \Sigma\subseteq \Sigma$ (as $b\in\Sigma$), so
  $e\in \ua b\cap D\subseteq (\ua b\cap \da \Sigma)\cap D\subseteq \Sigma\cap D$, which implies that ${x}_{m,n}\in\da e\subseteq \da (\Sigma\cap D)$. 
Thus, $D\cap (X\setminus\maxx X)\subseteq \da (D\cap\Sigma)$. Therefore, $D\subseteq \da (D\cap \Sigma)$. 
This shows that $D\cap \Sigma$ is cofinal in $D$,  so $D\cap \Sigma$ is a directed set such that $\bigvee D=\bigvee (D\cap \Sigma)$ (see  pp.61 in \cite{goubault}).
\medskip

(5) $L$ is an $\omega$-algebraic domain and $K(L)=L\setminus(\maxx(\Sigma)\cup\maxx(\Sigma^*)\cup\maxx X)$.

Let $B:=L\setminus(\maxx(\Sigma)\cup\maxx(\Sigma^*)\cup\maxx X)$.
Then, $B=\{x_{m,n}:m,n\in\mn\}\cup\{a\in\Sigma: \ell(a)\in\mn\}\cup \{a^*:a\in\Sigma, \ell(a)\in\mn\}$. Suppose $m\in\mn$ and $a=\langle{a_1,a_2,\cdots,a_m}\rangle\in\Sigma$.

\medskip

\emph{Step 1: } We prove $a^*\in K(L)$.

Suppose $D$ is a directed subset of $L$ such that ${a^*}\leq \bigvee D$. Then, $\bigvee D\in \ua {a^*}\subseteq \ws$, so
by Claim 3 of (3), $D\cap \Sigma^*$ is a directed subset of $L$ and
$a^*\leq \bigvee D=\bigvee (D\cap \Sigma^*)$. It follows that $D\cap \Sigma^*$ is a directed subset of $\Sigma^*$ and
$a^*\sqsubseteq \bigvee D=\bigvee (D\cap \Sigma^*)$. Note that $a^*$ is a finite sequence of $\Sigma^*$, so it is compact in $(\Sigma^*,\sqsubseteq^*)$, and then there exists $b^*\in D\cap \Sigma^*$ such that $a^*\sqsubseteq^*b^*$. This implies that $a\leq b^*$ in $L$. Therefore, $a^*\in K(L)$.

%
%


\medskip

\emph{Step 2: } We prove $a\in K(L)$.

Suppose $D$ is a directed subset of $L$ such that ${a}\leq \bigvee D$. Then, $\bigvee D\in \ua {a}\subseteq \Sigma\cup\ws$. We consider the following two cases:

(e1) $\bigvee D\in \Sigma^*$.

Then, by Claim 2 of (1), $a^*\leq \bigvee D$, and as we have proved in Step 1 that $a^*\in K(L)$, there exists $e\in D$ such that $a\leq a^*\leq e$, so $a\leq e$.

(e2) $\bigvee D\in\Sigma$.

Then, by Claim 2 of (3), $D\cap \Sigma$ is a directed subset of $L$ and
$a\leq \bigvee D=\bigvee (D\cap \Sigma)$. It follows that $D\cap \Sigma$ is a directed subset of $\Sigma$ and
$a\sqsubseteq \bigvee D=\bigvee (D\cap \Sigma)$. Note that $a$ is a finite sequence of $\Sigma$, so it is compact in $(\Sigma,\sqsubseteq)$, and then there exists $e\in D\cap \Sigma$ such that $a\sqsubseteq e$. This implies that $a\leq e$ in $L$.

By (e1) and (e2), we conclude that  $a\in K(L)$.

\medskip

\emph{Step 3: } $\forall k, m\in\mn$,  ${x}_{k,m}\in K(L)$.

 As a matter of fact, assume that  $D$ is a directed subset of $L$ such that ${x}_{k,m}\leq \bigvee D$.
	
	(i) Let $\bigvee D\in X$. Then,
	 $D\subseteq \da X=X$. Note that $x_{k,m}\leq \bigvee D$ in $L$ implies $x_{k,m}\leq_X \bigvee D$ in $X$, and it is clear that $x_{k,m}\in K(X)$. Then, there exists $e\in D$ such that $x_{k,m}\leq_X e$, and hence $x\leq e$ in $L$.
	
	(ii) Let $\bigvee D\in{\Sigma}\cup \Sigma^*$. Then, as ${x}_{k,m}\leq \bigvee D$, there exist $n_1,n_2,\cdots,n_{k}\in\mn$ such that $n_k\geq m$ and    ${x}_{k,m}\leq\langle{n_1,n_2,\cdots,n_{k}}\rangle\leq  \bigvee D$.
	As we have proved in in Step 2, $\langle{n_1,n_2,\cdots,n_k}\rangle\in K(L)$,  there exists $e\in D$ such that $\langle{n_1,n_2,\cdots,n_k}\rangle\leq e$.
	Thus, $x_{k,m}\leq e$.
	
Therefore, by above (i)--(ii), ${x}_{k,m}$ is compact.

From Steps 1--3, we deuce that  $B\subseteq K(L)$. Clearly, every element of $L$ can be represented as the  supremum of some directed subset of $B$, which means that $B$ is a base of $L$. Note that $K(L)$ is included in any base of $L$, so $K(L)\subseteq B$.
Therefore $B=K(L)$.
Also, one easily see that $|B|=|K(L)|=\omega$. Hence, $L$ is an $\omega$-algebraic domain.

\medskip

	\item [(6)] $\maxx L$ is not a $G_{\delta}$-set in $\Sigma L$.
	
	Suppose, on the contrary, that $\maxx L$ is a $G_{\delta}$-set.
	Then, there exists a countable family $\{U_k:k\in\mn\}\subseteq \sigma(L)$ such that
	$\maxx L=\bigcap_{k\in\mn}U_k$. Recall that $\maxx L=\maxx(\Sigma^*)\cup \maxx X$.
	
	\emph{Procedure 1: } Since ${x}_{1,\omega}=\bigvee_{n\in\mn}{x}_{1,n}\in \maxx X\subseteq \maxx L\subseteq U_1\in\sigma(L)$, there exists $n_{1}\in \mn$ such that ${x}_{1,n_{1}}\in U_1$. Note that ${\langle}n_{1}{\rangle}\geq {x}_{1,n_{1}}\in U_1=\ua U_{1}$, so $${\langle}n_{1}{\rangle}\in U_1.$$
	
	\emph{Procedure 2: } Since ${x}_{2,\omega}=\bigvee_{n\in\mn}{x}_{2,n}\in \maxx X\subseteq \maxx L\subseteq U_2\in\sigma(L)$, there exists $n_{2}\in \mn$ such that ${x}_{2,n_{2}}\in U_2$. Note that ${{\langle}n_{1},n_{2}{\rangle}}\geq {x}_{2,n_{2}}\in U_2=\ua U_{2}$, so ${{\langle}n_{1},n_{2}
	\rangle}\in U_2$.
Since ${{\langle}n_{1},n_{2}{\rangle}}\geq {\langle}n_{1}\rangle\in U_1=\ua U_1$, we have that  ${{\langle}n_{1},n_{2}{\rangle}}\in U_1$, which follows that
	$${{\langle}n_{1},n_{2}{\rangle}}\in U_1\cap U_2.$$
	
		\emph{Procedure 3: } Since ${x}_{3,\omega}=\bigvee_{n\in\mn}{x}_{3,n}\in \maxx X\subseteq \maxx L\subseteq U_3\in\sigma(L)$, there exists $n_{3}\in \mn$ such that ${x}_{3,n_{3}}\in U_3$. Note that ${{\langle}n_{1},n_{2},n_{3}{\rangle}}\geq {x}_{3,n_{3}}\in U_3=\ua U_{3}$, so ${{\langle}n_{1},n_{2},n_{3}{\rangle}}\in U_3$. Since ${{\langle}n_{1},n_{2},n_3{\rangle}}\geq {{\langle}n_{1},n_2{\rangle}}\in U_1\cap U_2=\ua (U_1\cap U_2)$, we have that ${{\langle}n_{1},n_{2},n_3{\rangle}}\in U_1\cap U_2$, which follows that
		$${{\langle}n_{1},n_{2},n_3{\rangle}}\in U_1\cap U_2\cap U_3.$$

Repeating the above procedures, we obtain a countably infinite sequence $a: = {\langle}n_{k}{\rangle}_{k\in\mn}\in \bigcap_{k\in\mn}U_k$. Note that $a\notin\maxx L$ (as $a<a^*$),
 which contradicts our assumption $\maxx L=\bigcap_{k\in\mn}U_k$. Therefore, $\maxx L$ is not a $G_{\delta}$-set in $\Sigma L$.
\end{example}

From the preceding example, we obtain our main result as follows, which gives a negative answer to Martin's problem.

\begin{corollary}
There exists an $\omega$-algebraic domain $P$ such that $\mbox{Max}(P)$ is not a $G_\delta$-set of the Scott space $\Sigma P$.
\end{corollary}

\vspace{0.7cm}

\noindent{\bf Acknowledgments: }
This work was supported by the National Natural Science Foundation of China (1210010153, 12071188),  Jiangsu Provincial
	Department of Education (21KJB110008), and International Postdoctoral Exchange Fellowship Program (No. PC2022013).

\vspace{0.7cm}

\noindent{\bf Declaration of interests: }
The authors declare that they have no known competing financial interests or personal relationships that could have appeared to influence the work reported in this paper.

\end{document}